 \documentclass[final,3p,times]{elsarticle}

\usepackage{amssymb}
\usepackage{amsmath}
\usepackage{empheq}
\usepackage{esint}

 \usepackage{lineno}

\usepackage{color}

\def \ve {\varepsilon}
\journal{Applied Mathematics Letters}
\newtheorem{theorem}{Theorem}
\newtheorem{proposition}{Proposition}
\newtheorem{definition}{Definition} 
\newtheorem{remark}{Remark}
\begin{document}
\begin{frontmatter}



\title{Phase-Field Model of Cell Motility: Traveling Waves and Sharp Interface Limit}


\author[Penn]{Leonid Berlyand}
\author[Penn]{Mykhailo Potomkin}
\author[VR]{Volodymyr Rybalko}

\address[Penn]{Department of Mathematics, The Pennsylvania State University, University Park, PA 16802, USA}

\address[VR]{Mathematical Division, B. Verkin Institute for Low Temperature, Physics and Engineering of National Academy of Sciences of Ukraine, 47 Lenin Ave., 61103, Kharkiv, Ukraine}

\begin{abstract}
    This letter  is concerned with   asymptotic analysis of a PDE  model for motility of a eukaryotic cell on a substrate. This model was introduced in \cite{ZieSwaAra11}, where it was shown numerically that it successfully reproduces experimentally observed phenomena of cell-motility such as a discontinuous onset of motion and shape oscillations. The model consists of a  parabolic PDE for a scalar phase-field function coupled with a vectorial parabolic PDE for the actin filament network (cytoskeleton).
    We formally derive
   the sharp interface limit (SIL), which describes the motion of the cell membrane  and show that it is  a volume preserving curvature driven motion with an additional nonlinear term due to adhesion to the substrate and protrusion by the cytoskeleton.
     In a 1D model problem we rigorously justify the SIL, and, using numerical simulations, observe some surprising features such as discontinuity of
    interface velocities and hysteresis.
    We show that  nontrivial
    traveling wave solutions appear when the key physical parameter exceeds a certain critical value
 and the potential in the equation for phase field function possesses certain asymmetry. 
    
\end{abstract}


\end{frontmatter}


\section{Introduction}
\label{section:intro}

An initially symmetric cell on a substrate may exhibit 
spontaneous breaking of symmetry or 
self-propagation
along the straight line maintaining the same shape over many times of its length \cite{KerPinAllBarMarMogThe08,BarLeeKerMogThe11}.  Understanding the  initiation 
of steady motion of a biological cell  as well  as  the mechanism of symmetry breaking is a 
fundamental issue in cell biology.

In \cite{ZieSwaAra11,ZieAra13} a phase-field model was proposed
to describe motility of a eukaryotic cell on a substrate.
We consider a  simplified version of that
model without myosin contraction ($\gamma=0$ in  \cite{ZieSwaAra11}), which consists of two coupled PDEs
\begin{eqnarray}
&&\frac{\partial \rho_\ve}{\partial t}=\Delta \rho_\ve
-\frac{1}{\ve^2}W^{\prime}(\rho_\ve)
-P_\ve\cdot \nabla \rho_\ve +\lambda_\ve(t),\quad x\in \Omega, \;t>0,
\label{eq1}\\
&&\frac{\partial P_\ve}{\partial t}=\ve\Delta P_\ve -\frac{1}{\ve}P_\ve
-\beta \nabla \rho_\ve
\label{eq2}
\end{eqnarray}
in a bounded domain $\Omega\subset \mathbb{R}^2$,
where the unknowns are the phase-field function $\rho_\ve$ and the vector field $P_\ve$ modeling average orientation of the actin network. 
System \eqref{eq1}-\eqref{eq2} is obtained by diffusive
scaling of equations from \cite{ZieSwaAra11} to study a sharp interface limit (SIL) of that
model under special scaling assumptions on the parameters. 
We introduce the
volume preservation constraint via the Lagrange multiplier
\begin{equation} \label{lagrange}
\lambda_\ve(t)=\frac{1}{|\Omega|}\int_\Omega\left(\frac{1}{\ve^2}W^\prime(\rho_\ve)
+ P_\ve\cdot \nabla \rho_\ve \right)\, dx
\end{equation}
in place of the volume constraint originally introduced in the potential \cite{ZieSwaAra11}.
The function $W^{\prime}(\rho)$ in \eqref{eq1} is the derivative of a double equal well potential  (e.g.,  $W(\rho)=\frac{1}{4}\rho^2(1-\rho)^2$).

The phase-field function $\rho_\ve$ takes values close to the wells of
the potential $1$ and $0$ for sufficiently small $\ve>0$
everywhere in $\Omega$ except for a thin transition layer. The
corresponding subdomains are interpreted  as the inside cell and
the outside cell regions, while the transition layer models the cell membrane. In \eqref{eq2}, $\beta>0$ 
is a fixed parameter
responsible for the creation of the 
field $P_\ve$ near the interface. 
 The boundary conditions $\partial_\nu \rho_\ve=0$ and   $P_\ve=0$
are imposed  on the boundary $\partial \Omega$.



We study system \eqref{eq1}-\eqref{eq2} in the sharp interface limit $\ve\to 0$. 
 Well known approaches in the study of sharp interface limits of phase field models such as 
 viscosity solutions techniques and the $\Gamma$- convergence method,  see, e.g.,  
\cite{LioKimSle04,Gol97,Ser10},
 are not readily applied to \eqref{eq1}-\eqref{eq2} 
 because of the coupling through the terms $P_\ve\cdot \nabla \rho_\ve$ and 
 $\nabla \rho_\ve$. The comparison principle, necessary  for  the viscosity solutions technique, does not apply for \eqref{eq1}-\eqref{eq2}.  
  Also this system is not a gradient flow for an energy functional which makes the $\Gamma$-convergence techniques  inapplicable. 
 Another analytical approach, based on 
 formal asymptotic expansions  was developed for different phase field models in \cite{Che94,MotSha95,CheHilLog10}. 
 Some ingredients of this approach are also used in the present study. 
 We also mention here an alternative approach to cell motility  based on numerical study of free boundary value problems  
 developed in \cite{KerPinAllBarMarMogThe08,RubJacMog05,BarLeeAllTheMog15,RecTru13,RecPutTru15}, and numerical studies of different phase field models of cell motility \cite{CamZhaLiLevRap13}.

 In this work  we first show that solutions of  \eqref{eq1}-\eqref{eq2}
 do not blow up on finite time intervals for sufficiently small $\ve$ by establishing  energy type and pointwise bounds, 
 next we formally derive a law of motion of the interface postulating a two-scale ansatz in the spirit of \cite{MotSha95}. 
 Then we prove the existence of nontrivial traveling waves in a one-dimensional version
 of  \eqref{eq1}-\eqref{eq2} in the case when the potential $W$ has certain assymmetry. This is done by an asymptotic reduction to a finite dimensional system for $V$ and $\lambda$, and applying the 
 Schauder fixed point theorem. Finally in a one-dimensional dynamical system we rigorously 
 prove that the interface velocity satisfies a simple nonlinear equation and demonstrate existence of a hysteresis loop in the system by numerical simulations.    
%
%

\section{Existence of Solutions and Sharp Interface Limit in 2D Model}\label{sec_sil_2D}

The first result of this work demonstrates that for sufficiently
small $\ve>0$ a unique solution $\rho_\ve$, $P_\ve$ of
\eqref{eq1}-\eqref{eq2} exists and $\rho_\ve$ maintains the 
structure of a sharp interface between two phases $0$ and $1$,
provided that initial data are well prepared. To formulate this result
we introduce the following auxiliary (energy-type) functionals:
\begin{equation}
\label{energy}
\begin{array}{l}
E_\ve(t):=\frac{\ve}{2} \int_\Omega |\nabla \rho_\ve(x,t)|^2dx+\frac{1}{\ve}\int_\Omega W( \rho_\ve(x,t) )dx,\\ \\
F_{\ve}(t):= \int_\Omega \Bigl(| P_\ve(x,t)|^2+|P_\ve(x,t)|^4\Bigr)dx.
\end{array}
\end{equation}

\begin{theorem} \label{wp_theorem}  Assume that the system \eqref{eq1}-\eqref{eq2} is supplied with initial 
	data that satisfy $-\ve^{1/4}< \rho_\ve (x,0)<1+\ve^{1/4}$, and
    \begin{equation}
    E_\ve(0)+F_\ve(0)\leq C_1.
    \label{IniEnBound}
    \end{equation}
Then for any $T>0$ there exists  a solution $\rho_\ve$, $P_\ve$ of \eqref{eq1}-\eqref{eq2}  on the time interval $(0,T)$ when  $\ve>0$ is sufficiently small, $\ve<\ve_0(T)$.
    Moreover, $-\ve^{1/4}\leq \rho_\ve (x,t)\leq 1+\ve^{1/4}$ and
    \begin{equation}\label{noblowup1}
    \ve \int_0^T\int_\Omega \Bigl(\frac{\partial\rho_\ve}{\partial t} \Bigr)^2dxdt\leq C_2,
    \quad
    E_\ve(t)+F_\ve(t)\leq C_2\quad \forall t\in(0,T),
    \end{equation}
    where $C_2$ is independent of $t$ and $\ve$.
\end{theorem}

This theorem shows that there is no blow up of the solution on the given time interval $(0,T)$,
also it proves that if the initial data have sharp interface structure, this sharp interface structure is 
preserved by the solution on the whole time interval $(0,T)$. 
The claim of Theorem \ref{wp_theorem}  is nontrivial due
to the presence of the quadratic term
$P_{\ve}\cdot \nabla\rho_{\ve}$ in \eqref{eq1} which, in general, could lead to a finite time blow up.
The main idea behind the existence proof is to find and
utilize a bound for $\rho_\ve$ in $L^\infty((0,T)\times \Omega)$,
which is obtained by combining the maximum principle and energy estimates.


Next we study the SIL $\ve\to 0$  for the system (\ref{eq1})-(\ref{eq2}). 
 We seek solutions in the
form of ansatz (locally in a neighborhood of the interface)
\begin{equation}
\label{ANSATZ}
\rho_\ve=\theta_0(d/\ve)+\ve
\theta_1(d/\ve,S) +\dots,\quad
P_\ve=\nu\Psi_0(d/\ve,
S)+\dots,
\end{equation}
where $d=d(x,t)$ is the (signed) distance to a unknown evolving
interface curve $\Gamma(t)$, $S=s(p(x,t),t)$ with
$p(x,t)$ being the projection of
$x$ on $\Gamma(t)$ and $s(\xi,t)$
being a parametrization of
$\Gamma(t)$, $\nu=\nu(p(x,t),t)$ is the inward pointing normal to
$\Gamma(t)$ at $p(x,t)\in\Gamma(t)$. The key choice 
here is the interface curve $\Gamma(t)$ that allows for 
appropriate estimates.
We substitute this ansatz in
(\ref{eq1}) to find, after collecting terms (formally) of the order $\ve^{-2}$, that $\theta_0$ satisfies $\theta_0^{\prime\prime}=W^\prime(\theta_0)$.
It is known that there exists a unique (up to a translation) solution (standing wave)
$\theta_0(z)$ which tends
to $0$ or $1$ when $z\to-\infty$ or $z\to+\infty$. For the potential
$W(\rho)=\frac14\rho^2(\rho-1)^2$ the function $\theta_0$ is explicitly given by
$\theta_0(z)=\frac{1}{2}\left(1+\tanh \frac{z}{2\sqrt{2}}\right)$.  Then substitute \eqref{ANSATZ} in (\ref{eq2}) and consider the leading
(of the order $\ve^{-1}$) term.  Denoting by $V(x,t)$ the (inward) normal velocity 
of the curve $\Gamma(t)$ at $x\in \Gamma(t)$
we  obtain that the scalar function
$\Psi_0(z)$ solves
\begin{equation}\label{eq_for_Psi}
-\frac{\partial^2\Psi_0}{\partial z^2}-V\frac{\partial \Psi_0}{\partial z}+\Psi_0+
\beta\theta_0^\prime(z)=0.
\end{equation}
Finally, assuming that the leading term of the expansion of $\lambda_\ve$ is of the order ${\ve}^{-1}$, $\lambda_\ve =\lambda(t)/\ve+\dots$,  and collecting terms of the order $\ve^{-1}$ in \eqref{eq2} we are led to the following equation
$$
-\frac{\partial^2 \theta_1}{\partial
z^2}+W^{\prime\prime}(\theta_0)\theta_1=(V-\kappa)\frac{\partial
\theta_0}{\partial z}-\Psi_0\frac{\partial \theta_0}{\partial
z}+\lambda(t),
$$
where $\kappa$ denotes the curvature of $\Gamma(t)$.
The solvability condition for this equation
(orthogonality to the eigenfunction  $\theta_0^\prime$ of the linearized Allen-Cahn equation) yields
the desired sharp interface equation
\begin{equation}
\label{SharpInterEq}
V(x,t)= \kappa(x,t) +
\frac{1}{c_0}\Phi_\beta(V(x,t))-\lambda(t),\quad x\in\Gamma(t),
\end{equation}
where $c_0=\displaystyle\int\left(\theta_0^\prime\right)^2dz$, and $\Phi_\beta(V)$ is given by
\begin{equation} \label{def_of_Phi}
\Phi_\beta(V)=\int\limits_{\mathbb R} \Psi_0\left(\theta_0^\prime(z)\right)^2dz.
\end{equation}
From the volume preservation condition
$\int_{\Gamma(t)}Vds=0$ it follows that
$\lambda(t)=\frac{1}{c_0}\fint_{\Gamma(t)}(c_0\kappa+\Phi_\beta(V))ds$.

The above formal derivation of the sharp interface limit is rigorously
justified in 1D (see Theorem \ref{1Dinterface} below) because of significant technical difficulties due to the curvature in 2D.   
Solvability of \eqref{SharpInterEq} was shown
in \cite{MizBerRybZha15} for $\beta$ less than some
critical value, moreover \eqref{SharpInterEq} was proved to enjoy a parabolic 
regularization feature.  However  for large $\beta$, the equation \eqref{SharpInterEq} might have multiple solutions.  
To obtain a selection criterion 
and 
elucidate the role of the parameter $\beta$ in the cell interface motion we consider a 1D 
model of the cell-motility in the next sections.


\section{Traveling wave solutions in 1D}\label{sec_tw_1D}

In this section we 
show that solutions of system \eqref{eq1}-\eqref{eq2}  exhibit
significant qualitative changes when the parameter $\beta$
increases and the potential $W(\rho)$ has certain asymmetry, e.g.  $W(\rho)=\frac14\rho^2(\rho-1)^2(1+\rho^2)$. Here we look for traveling wave solutions in 1D model,
considering (\ref{eq1})-(\ref{eq2}) with $\Omega=\mathbb{R}^1$. In other words we
are interested in nontrivial spatially localized solutions of \eqref{eq1}-\eqref{eq2} of the form
$\rho_{\ve} =\rho_{\ve}(x-Vt)$, $P_{\ve}=P_{\ve}(x-Vt)$. This leads to the stationary
equations with unknown (constant) velocity $V$ and constant $\lambda$:
\begin{eqnarray}
0&=&\partial^2_x \rho_{\ve} +V\partial_x\rho_{\ve}-\frac{W'(\rho_{\ve})}{\ve^2}-P_{\ve}\partial_x\rho_{\ve} +\dfrac{\lambda}{\ve}\label{tw_rho},\\
0&=&\ve \partial_x^2P_{\ve}+V\partial_xP_{\ve} -\frac{1}{\ve} P_{\ve} -\beta\partial_x\rho_{\ve}. \label{tw_P}
\end{eqnarray}

We are interested in solutions of \eqref{tw_rho}-\eqref{tw_P} that are essentially localized on the interval $(-a,a)$, for a  given $a>0$. 
We look for such solutions for sufficiently small
$\ve>0$ with the phase field function $\rho_\ve$ of the form
\begin{equation}
	\label{repr}
\rho_{\ve} =\theta_0((x+a)/\ve)\theta_0((a-x)/\ve)+\ve\psi_\ve+\ve u_{\ve},
\end{equation}
where constant $\psi_\ve$ is the smallest solution of 
$W^\prime(\ve\psi)=\ve\lambda$ and $u_{\ve}$ is the
new unknown function vanishing at $\pm\infty$. Observe that the first term 
$\theta_0((x+a)/\ve)\theta_0((a-x)/\ve)$ has "$\Pi$" shape and becomes the characteristic 
function of the interval $(-a,a)$ in the limit $\ve\to 0$.

\begin{proposition} \label{prop_1}
	For any real $\beta\geq 0$ and sufficiently small $\ve$ there exists a localized standing wave  solution (with $V=0$) of \eqref{tw_rho}-\eqref{tw_P} .  It is localized in the sense that  the
	representation \eqref{repr} holds with $u_{\ve}\in L^2(\mathbb{R})\cap L^\infty(\mathbb{R})$ and $\|u_{\ve}\|_{L^\infty}\leq C$.  
\end{proposition}
Proposition \ref{prop_1}  justifies expected  existence of standing wave solutions  (immobilized cells)  in the class of functions with the symmetry 
$\rho(-x)=\rho(x)$ and $P(-x)=-P(x)$, so that the polarization field on the front and back has the same magnitude but is oriented in opposite directions. This field, loosely speaking, is trying to push front and back in opposite directions with the same velocities, thus, cell does not move. 
Indeed, the relation between $P_\ve$ and $V$ can be obtained from
the second equation in \eqref{ANSATZ}, \eqref{def_of_Phi} and \eqref{tw_front_and_back}.

We show, however, that not all localized  solutions of \eqref{tw_rho}-\eqref{tw_P}
are necessarily standing waves. 
Assuming that there exists a traveling wave solution with a nonzero velocity, e.g. $V>0$, and passing to 
the sharp interface limit $\ve\to 0$ in \eqref{tw_rho}-\eqref{tw_P} 
at the back  and front  transition layers ($x=\pm a$ in \eqref{repr}) 
we formally obtain two relations for the velocity $V$ and the constant $\lambda$  
\begin{equation}\label{tw_front_and_back}
c_0V=\Phi_\beta(V)-\lambda, \text{  and  } -c_0V=\Phi_\beta(-V)-\lambda. 
\end{equation} 
Then eliminating  $\lambda$ we obtain the equation for the velocity $V$:
\begin{equation}\label{eq_for_nonzero_tw}
 2c_0V =\Phi_\beta(V)-\Phi_\beta(-V).
\end{equation} 
This equation always has one root $V=0$ which corresponds to the standing wave solution whose existence for system \eqref{tw_rho}-\eqref{tw_P} 
is established in Proposition \ref{prop_1}.  Two more roots, say $V_0$, and $-V_0$ appear for 
sufficiently large $\beta>0$ in the case when $\Phi_\beta(V)>\Phi_\beta(-V)$ for $V>0$, thanks to the fact that
$\Phi_{\beta}$ is proportional to $\beta$ (note that if $W(\rho)=\frac{1}{4}\rho^2(\rho-1)^2$
then $\Phi_\beta$ is an even function, so the RHS of \eqref{eq_for_nonzero_tw} vanishes for arbitrary $\beta$ and thus $V$ is necessarily $0$). 
 This heuristic argument can be made rigorous by proving the following:
\begin{theorem} \label{theorem_2}
	Let $W(\rho)$ and $\beta$ be such that \eqref{eq_for_nonzero_tw}
	has a root $V=V_0>0$ and $\Phi_\beta^\prime(V_0)+\Phi_\beta^\prime(-V_0)\not=2c_0$
	(nondegenerate root).
	 Then for sufficiently small $\ve>0$ 
	there exists a localized solution of \eqref{tw_rho}-\eqref{tw_P} with $V=V_\ve\neq 0$, 
	moreover
	$V_{\ve}\to V_0\neq 0$ as $\ve \to 0$ (as above localized solution means that representation \eqref{repr} holds with $u_{\ve}\in L^2(\mathbb{R})\cap L^\infty(\mathbb{R})$ and $\|u_{\ve}\|_{L^\infty}\leq C$). 
\end{theorem}

\noindent{\bf Remark.} In Theorem 2, it is crucial that  \eqref{eq_for_nonzero_tw} has a non-zero solution $V_0$ which is impossible for the symmetric potential $W(\rho)=\frac{1}{4}\rho^2(\rho-1)^2$, but does hold for an asymmetric potential, e.g., $W(\rho)=\frac{1}{4}\rho^2(\rho-1)^2(1+\rho^2)$. 
In the case of smaller diffusion in equation \eqref{tw_P} one can prove that $\int_0^1 W''(\rho)dW^{3/2}(\rho)>0$ is a sufficient condition for existence of $V_0\neq 0$. We conjecture that this remains true for \eqref{tw_rho}-\eqref{tw_P}.

Theorem \ref{theorem_2} guarantees existence of non-trivial traveling waves that describe steady motion without external stimuli. 
Thus our analysis of \eqref{tw_rho}-\eqref{tw_P} is consistent with 
experimental observations of motility on keratocyte cells \cite{KerPinAllBarMarMogThe08}.

The proof of Theorem \ref{theorem_2} is carried out in two steps. 
In the first step we use \eqref{repr} to rewrite \eqref{tw_rho}-\eqref{tw_P} as a single equation of the form $\mathcal{A}_{\ve}u_{\ve}+\ve B_\ve(V,\lambda)+\ve^2 C_\ve(u_{\ve},V,\lambda)=0$, where $\mathcal{A}_\ve u:=\ve^2\partial_x^2u-
W^{\prime\prime}(\theta_0((x+a)/\ve)\theta_0((a-x)/\ve))u$ is the Allen-Cahn operator linearized around the first term in \eqref{repr}. We rewrite this equation as a fixed point problem $u_{\ve}=-\ve\mathcal{A}_{\ve}^{-1}( B_\ve(V,\lambda)+\ve C_\ve(u_{\ve},V,\lambda))$. The operator $\mathcal{A}_{\ve}$ has zero eigenvalue of multiplicity two (up to a proper $o(\ve^2)$ perturbation).
This leads to solvability conditions 
which to the leading term coincide with \eqref{tw_front_and_back}. In the second step we apply the Schauder fixed point theorem to establish existence of solutions of \eqref{tw_rho}-\eqref{tw_P}. 

%


\section{Sharp interface limit in a  1D model problem and hysteresis}
\label{sec_sil_1D}

This section is devoted to the asymptotic analysis as $\ve \to 0$ of the following 1D problem
\begin{eqnarray}
&&
\frac{\partial \rho_{\varepsilon}}{\partial t}=\partial^2_{x}\rho_{\varepsilon}-\frac{W'(\rho_{\varepsilon})}{\varepsilon^2}-P_{\varepsilon}\partial_x\rho_{\varepsilon}+\frac{F(t)}{\varepsilon},  \label{1D_rho}
\\
&&\frac{\partial P_{\varepsilon}}{\partial t}=\varepsilon \partial_{x}^2P_{\varepsilon}-\frac{1}{\varepsilon}P_{\varepsilon}-\beta \partial_{x}\rho_{\varepsilon},
\label{1D_P}
\end{eqnarray}
$x\in \mathbb{R}^1$, $t>0$, for a given function $F:(0,+\infty)\to \mathbb{R}^1$. This is a model 
problem to develop rigorous mathematical tools for  \eqref{eq1}-\eqref{eq2}, and it describes a normal cross-section of the transition layer (interface) between $0$ and $1$ phases. The variable $x\in \mathbb{R}$ corresponds  to the re-scaled signed distance $d$ (see Section \ref{sec_sil_2D}). The function 
$F(t)$ models forces due to the curvature 
of the interface and the mass preservation constraint $\lambda_\ve$, and for technical simplicity $F(t)$ is chosen to be independent of $x$.


Similar to Section \ref{sec_tw_1D}, we seek the solution of \eqref{1D_rho}-\eqref{1D_P} in the form 
\begin{equation}\label{eq_form}
\rho_{\ve}(x,t)=\theta_0(y)+\ve\psi_{\ve}(y,t)+\ve u_{\ve}(y,t), \;\;y=\frac{x-x_{\ve}(t)}{\ve},
\end{equation}

\vspace{-0.05 in}

\noindent where $\theta_0$ and $\psi_{\ve}$ are known functions, and $u_{\ve}$ is a new unknown function. Function $\psi_{\ve}(y,t)$ is defined by
\begin{equation}\nonumber
\psi_{\ve}(y,t)=\psi^-_{\ve}(t)+\theta_0(y)(\psi^+_{\ve}(t)-\psi^-_{\ve}(t)), \quad\text{ where}\quad  \partial_t (\ve \psi^{\pm}_{\ve})=-\frac{W'((1\pm 1)/2+\ve\psi_{\ve}^{\pm})}{\ve^2}+\frac{F(t)}{\ve},\; \psi^{\pm}_{\ve}(0)=0.
\end{equation}
 Existence of the 
 $x_{\ve}(t)$ (describing the location of the interface) together with estimates on $u_{\ve}$ uniform in $\ve$ and $t$ are established in the following

\begin{theorem}\label{ansatz_existence}
	Let $\rho_{\ve},P_{\ve}$ be a solution of Problem \eqref{1D_rho}-\eqref{1D_P} with initial data for $\rho_\ve$ and $P_{\ve}$ satisfying 
	"well-prepared"  initial conditions:
	\begin{equation}\label{ic}
	\rho_{\ve}(x,0)=\theta_0\left(x/{\ve}\right)+\ve v_{\ve}\left({x}/{\ve}\right),
	\end{equation}
	where $\|v_{\ve}\|^2_{L^2}=\int_{\mathbb R}|v_{\ve}(y)|^2 dy<C$, $\|v_\ve\|_{L^{\infty}(\mathbb R)}\leq C/\ve$, and $P_{\ve}(x,0)=p_\ve (\frac{x}{\ve})$ such that
	\begin{equation}
	\|p_\ve\|_{L^2(\mathbb R)}+ \| p_\ve\|_{L^{\infty}(\mathbb R)}+\| \partial_y p_\ve\|_{L^{\infty}(\mathbb R)}<C. 
	\end{equation}
	Then there exists 
	$x_{\ve}(t)$ such that expansion \eqref{eq_form}  holds with $\|u_{\ve}(\cdot,t)\|_{L^2(\mathbb R)}~<~C$ for $t\in [0,T]$ and $\int_{\mathbb R} u_{\ve}\theta_0' dy =0$. Moreover, 
	assuming that $\int_{\mathbb R} v_{\ve}\theta_0' dy =0$ , 
	the interface velocity $V_\ve=\dot x_\ve(t)$ is determined by the following system:
	 \begin{empheq}[left=\empheqlbrace]{align}
	(c_0+\ve \tilde{\mathcal{O}}_\ve (t))V_\ve(t)&=\int (\theta_0')^2 \Psi_{\ve}dy-F(t)+\ve\mathcal{O}_{\ve}(t),\label{eq_for_reduced_V}\\
	\ve\frac{\partial \Psi_{\ve}}{\partial t}&= \frac {\partial^2 \Psi_{\ve}}{\partial y^2}+V_{\ve}(t)\frac{\partial \Psi_{\ve}}{\partial y}-\Psi_{\ve}-\beta {\theta_0'}(y),\label{eq_for_reduced_A}
	\end{empheq}
	where ${\tilde{\mathcal{O}}}_{\ve}(t)$ and   $\mathcal{O}_{\ve}(t)$ are bounded in $L^{\infty}(0,T)$. 
\end{theorem}

The reduced system \eqref{eq_for_reduced_V}-\eqref{eq_for_reduced_A} can be further simplified by
taking the limit $\ve\to 0$. Formal passing to the limit in \eqref{eq_for_reduced_A} leads to
equation \eqref{eq_for_Psi} whose unique solution depends on the parameter $V$. 
Substituting this solution into \eqref{eq_for_reduced_V} in place of $\Psi_\ve$ we obtain the equation 
\begin{equation}\label{1D_sil}
c_0V_0(t)=\Phi_\beta(V_0(t))-F(t)
\end{equation}  
for the limiting velocity $V_0=\lim_{\ve\to 0} V_\ve$. However, in general, equation 
\eqref{1D_sil} is not uniquely solvable. The plot of the function $c_0V-\Phi_\beta(V)$ for sufficiently large $\beta$ is depicted 
on the Figure 1, where one sees that  \eqref{1D_sil} has two or three solutions when $F\in [F_{\rm min},F_{\rm max}]$. In order to justify \eqref{1D_sil}  and select a correct solution we reduce system \eqref{eq_for_reduced_V}-\eqref{eq_for_reduced_A} to a single nonlinear equation substituting expression for $V_\ve$ from  \eqref{eq_for_reduced_V} into \eqref{eq_for_reduced_A}. 
Then rescaling time and neglecting terms of the order $\ve$ we arrive at the equation
$\partial_t U=\partial_y^2 U+\frac{1}{c_0}(
 \int (\theta^\prime_0)^2Udy-F)\partial_yU-
 U-\beta \theta_0^\prime$ whose long time behavior has to be analyzed in order to obtain the limit of \eqref{eq_for_reduced_V}-\eqref{eq_for_reduced_A}  as $\ve\to 0$. This is done by spectral analysis of the linearized operator $\mathcal{A}_V U=\partial_y^2 U+V\partial_y U-
 U-\frac{1}{c_0}\partial_y \Psi_0
 \int  (\theta^\prime_0(z))^2U(z)dz$ about steady states $\Psi_0$ of the above nonlinear equation, 
 where $\Psi_0$ are obtained by finding  roots $V$ of the ordinary equation $c_0V=\Phi_\beta(V)-F$ 
 and then solving the PDE 
  \eqref{eq_for_Psi}.

\begin{definition} \label{def_stable} 
	Define the set of stable velocities $\mathcal{S}$ by
$\mathcal{S}=\{V \in \mathbb{R};\ \sigma(\mathcal{A}_V)\subset \{ \lambda\in \mathbb{C};{\rm Re}\lambda <0\}\}$, where $\sigma(\mathcal{A}_V)$ denotes the spectrum of the operator $\mathcal{A}_V$ (note that $\mathcal{S}$ is an open set).
\end{definition}

\begin{theorem}
	\label{1Dinterface}
	 Let $F(t)$ be a continuous function and assume that  $V_0\in \mathcal{S}$ solves $c_0V_0=\Phi_\beta(V_0)-F(0)$. Assume also that $\|p_\ve-\Psi_0\|_{L^2} \leq\delta$, where $\Psi_0$ is the solution of   \eqref{eq_for_Psi} with $V=V_0$ and $\delta>0$ is some small number depending 
on $V_0$ but independent of $\ve$.  Then $V_\ve(t)=\dot x_\ve(t)$ defined in Theorem
\ref{ansatz_existence} converges to the continuous solution of  the equation 
$c_0V(t)=\Phi_\beta(V(t))-F(t)$ with $V(0)=V_0$ on every finite time interval $[0,T]$ where such a 
solution exists and $V(t)\in\mathcal{S}$ $\forall t\in [0,T]$.
\end{theorem}

We conjecture that stability of velocities is related to monotonicity intervals of the function 
$c_0V-\Phi_\beta(V)$. This conjecture is supported by the following result.

\begin{proposition} If $c_0\leq\Phi^\prime_\beta(V)$, then $V$ is not a stable velocity. 
\end{proposition}

In general $\Phi^\prime_\beta(0)$ is nonzero if the potential $W(\rho)$ 
is asymmetric. In particular,  for $W(\rho)=\frac{1}{4}\rho^2(1-\rho)^2(1+\rho^2)$ we have 
$c_0<\Phi^\prime_\beta(0)$ when $\beta>\beta_{critical}>0$, therefore
zero velocity is not stable in this case. For 2D problem this would imply instability of initial circular shape leading to a spontaneous breaking of symmetry observed in experiments. 

\begin{remark}
	In the particular case $W(\rho)=\frac{1}{4}\rho^2(\rho-1)^2$ we prove that $(-\infty,\sqrt{2})\cap \left\{V;\; c_0>\Phi_\beta^\prime(V) \right\}\subset \mathcal{S}$. We also establish  $\mathcal{S}=\left\{V;\; c_0>\Phi_\beta^\prime(V) \right\}$ via verifying numerically a technical  inequality. 
	
	\end{remark}

%

While Theorem \ref{1Dinterface} describes local in time continuous evolution of
the interface velocity according to the law $c_0V=\Phi_\beta(V)-F(t)$
until $V$ leaves the set of stable velocities $\mathcal{S}$, we
conjecture that this law remains valid even after the time when
the solution $V$ reaches an endpoint of a connected component of
$\mathcal{S}$. Consider a particular example 
of
$\beta=150$, the corresponding plot of the function $c_0V-\Phi_\beta(V)$
is depicted on Fig. 1. \textcolor{black}
{Choose $F(t)$ given by
$F(t)=F_{\uparrow}(t):=-2.25+1.25t$ for $t\in [0,1]$ and
$F(t)=F_{\downarrow}(t):=F_{\uparrow}(2-t)$ for $t\in (1,2]$. Starting with well prepared
initial data we expect that the interface velocity $V$ increases
with $F(t)$ until it reaches $V_{\rm max}$ then it jumps to
another branch and continues to vary in $(V_{\rm min},+\infty)$
till the moment when it decreases to $V_{\rm min}$ and experiences
one more jump, then it varies in $(-\infty,V_{\rm max})$ to return to
the initial velocity at $t=2$ see Fig. 1, left.} Thus we conjecture that system has
a hysteresis loop, this conjecture is verified by numerical
simulations for the sharp interface limit \eqref{1D_sil} 
as well as the
original system \eqref{1D_rho}-\eqref{1D_P} for small $\ve$. The results
of the latter simulations with $\ve=0.01$ are depicted on Fig. 1, right.

\begin{figure}[t]
	\begin{center}
	\includegraphics[width=0.3\textwidth]{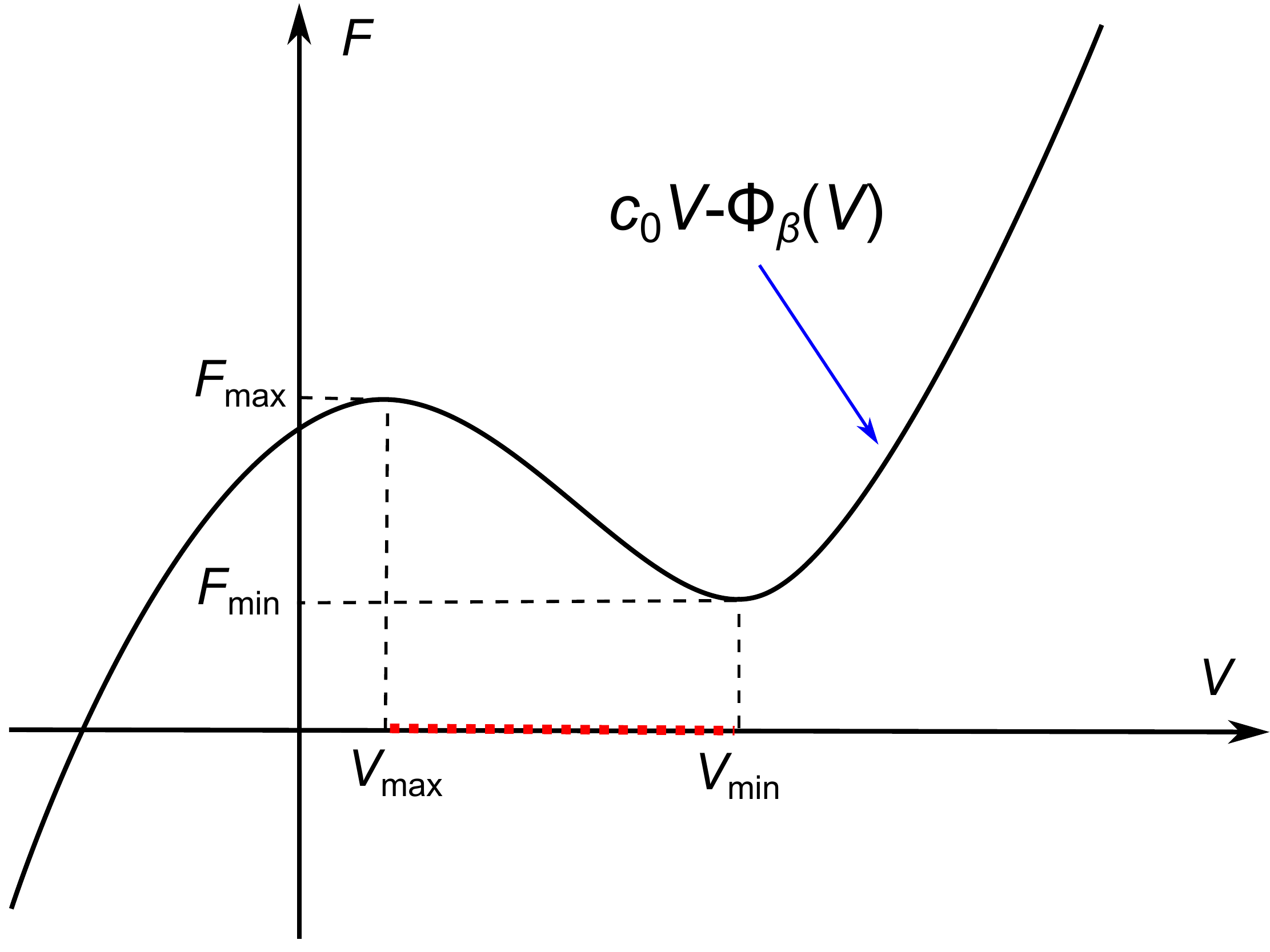}
	\includegraphics[width=0.34\textwidth]{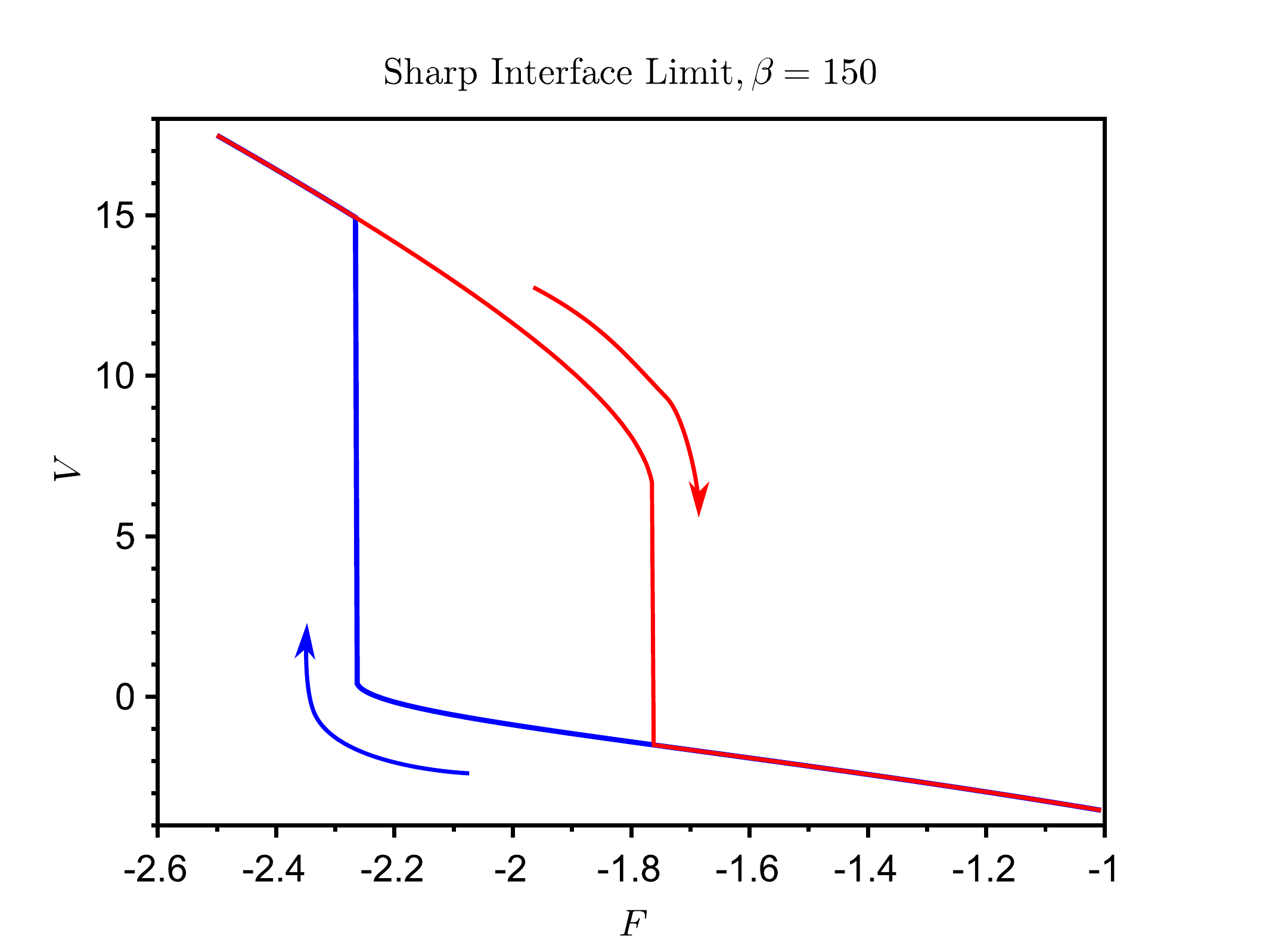}
	\includegraphics[width=0.34\textwidth]{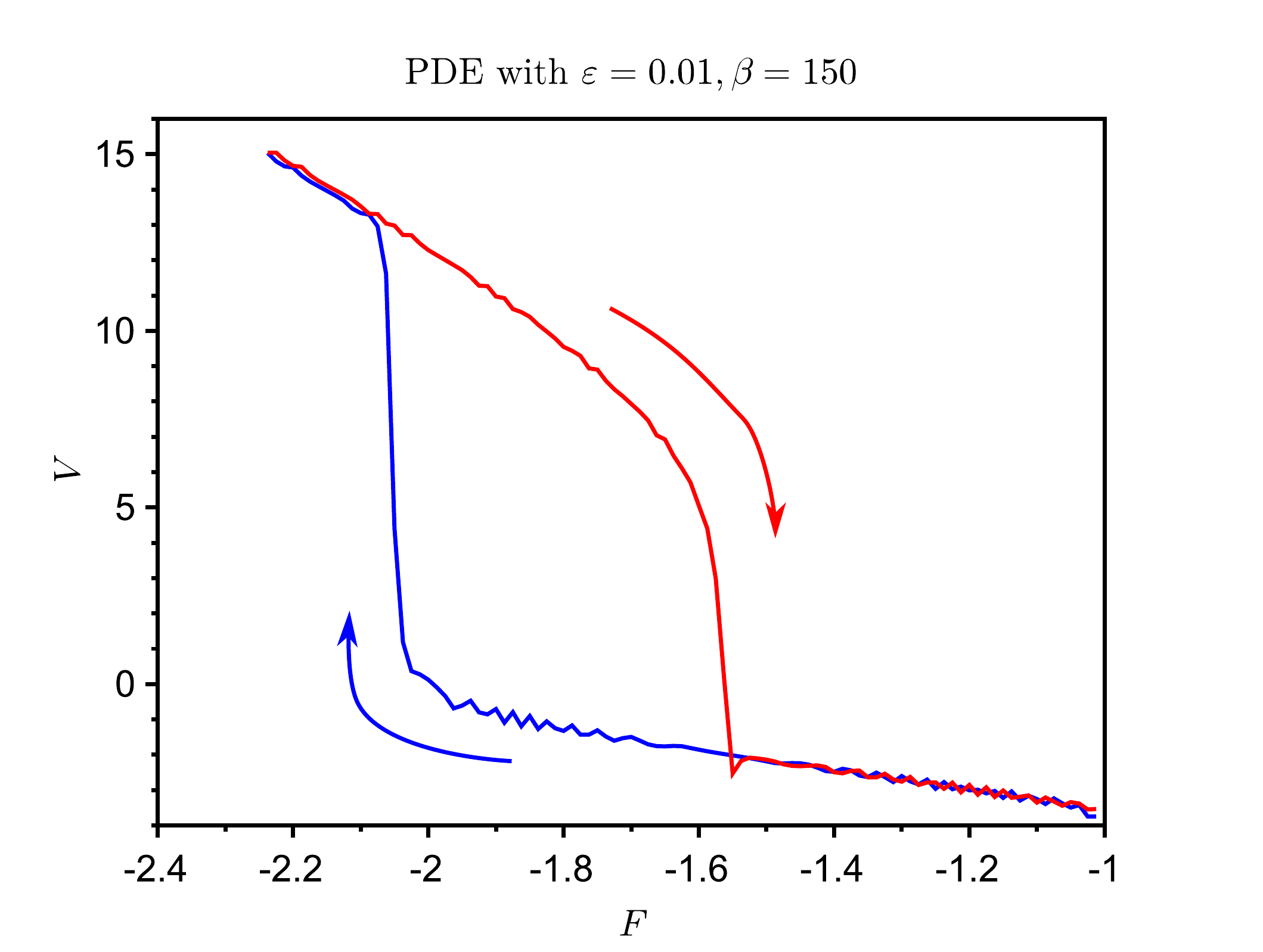}
	\caption{Hysteresis loop  in the problem of cell motility. (Left) The sketch of the plot for $c_0V- \Phi_\beta(V)$;
		(Center,Right) Simulations of $V=V(F)$, 
		(Center): solution of \eqref{SharpInterEq} 
		 (Right): solution of PDE system \eqref{1D_rho}-\eqref{1D_P}. On both figures (Center) and (Right) arrows show in what direction the system $(V(t),F(t))$ evolves as time $t$ grows; blue curve is for $F_{\downarrow}(t)$, red curve is for $F_{\uparrow}(t)$.}
	\label{fig:hysteresis}
	\end{center}
	\vspace{-0.15 in}
\end{figure}

\section*{Acknowledgments}

This work of LB and VR  was partially  supported by NSF grants  
DMS-1106666 and DMS-1405769. The work of MP was partially supported by the NSF grant DMS-1106666.

\section*{References}
\bibliographystyle{elsarticle-num}
\bibliography{cell}
\end{document}